\documentstyle[12pt]{amsart}
\title[The geometric minimal models of analytic spaces]{The geometric
minimal models\\
of
 analytic spaces}
\author{Shihoko Ishii}
\author{Pierre Milman}
\address{Shihoko Ishii: Department of Mathematics, Tokyo Institute of
Technology, Oh-Okayama, Meguro, Tokyo, Japan}
\address{Pierre Milman: Department of Mathematics, University of
Toronto, Toronto, Ontario, Canada}

\newcommand{\bC}{{\Bbb C}}
\newcommand{\bP}{{\Bbb P}}
\newcommand{\bZ}{{\Bbb Z}}
\newcommand{\bQ}{{\Bbb Q}}

\newtheorem{thm}{Theorem}[section]

\newtheorem{lem}[thm]{Lemma}
\newtheorem{cor}[thm]{Corollary}
\newtheorem{prop}[thm]{Proposition}

\theoremstyle{definition}
\newtheorem{defn}[thm]{Definition}

\newtheorem{say}[thm]{}
\newtheorem{exmp}[thm]{Example}

\newtheorem{rem}[thm]{Remark}

\theoremstyle{remark}

\begin{document}
\maketitle
\begin{abstract}
  This paper shows that an analytic space $X$
  has
 a unique maximal model  through which every proper surjective morphism from
 a non-singular analytic
  space to $X$ factors.
  This is called the {\sl geometric minimal model} of $X$ and
  characterized by the contraction property of rational curves.
  Some other properties such as functoriality, the direct product
  property and the
  quotient property of the geometric minimal model are also
  studied here.
  The relation of the geometric minimal model
  with  Mori's minimal model  is discussed.
\end{abstract}

\section{Introduction}
\noindent
  It is well known that a two dimensional singularity
  $(X,x)$ has a good birational model called the minimal resolution
  $f:Y \to X$,
  that is :

  (1) all resolutions of $(X, x)$ factor through $f$
  and

  (2) $Y$ itself is non-singular.

  For the higher dimensional case,
  such a resolution no longer exists.
  Then,  seeking a good birational model of a higher dimensional
  singularity
   in the Minimal Model Conjecture,
   people modify the conditions (1) and (2) into a condition (3)
   which is equivalent to (1) and  (2) in the two-dimensional case.

  (3) $Y$ has $\bQ$-factorial terminal singularities and $K_Y$ is $f$-nef.

  As is known, a model satisfying (3) exists for a three-dimensional
  singularity (\cite{mori}), but the existence for the higher
 dimensional case is still  a conjecture.

  On the other hand, the minimal resolution $f:Y\to X$ of a two-dimensional
  singularity is also characterized as follows:

 (1) all resolutions of $(X, x)$ factor through $f$
  and

 $(2')$ $f:Y\to X$ is maximal among those satisfying (1).

  This paper shows that such a model also uniquely exists for the
  higher dimensional
  case.

  We also prove that a model $f:Y\to X$ satisfying the following
  conditions \( (1') \), \( (2'') \)  uniquely exists for any
  dimension:

  $(1')$ all proper surjective morphisms  $\tilde X \to X$
  from non-singular $\tilde X$
   factor through $f$ and

  $(2'')$  $f:Y\to X$ is maximal among those satisfying $(1')$.

  This latter model is called the {\sl geometric minimal model} of $X$ and
is proved to
coincide
  with the model satisfying (1) and $(2')$ if $\dim X\neq 2$.

  The geometric minimal model is characterized by the contraction of rational
  curves:
\noindent
  let $\tilde X \to X$ be a proper surjective morphism from non-singular
$\tilde X$,
  then every rational curve on $\tilde X$
  which is mapped to a point on X is mapped to a point on the geometric minimal
  model $Y$;
  conversely, the maximal model $Y\to X$ satisfying the following property
  coincides with the geometric minimal model:
  there is a proper bimeromorphic morphism $Z\to Y$ such that every
rational curve
  on $Z$ which is mapped to a point on $X$ is mapped to a point on $Y$.
  By this, we can see that the geometric minimal model of a toric variety
  is the identity and the geometric minimal model of a cone over an Abelian
  variety is the corresponding line bundle over the Abelian variety.

  In \S 2, we prove the existence of a maximal model
  through which all resolutions factor.
  The proof is done  under a more general
  setting, and the existence of the geometric minimal model
  is also proved.

  In \S 3,  we prove the characterization of the geometric minimal model
  by the contraction of rational curves.
  The coincidence of the geometric minimal model and the model
  satisfying
  (1) and $(2')$ for $\dim \neq 2$ is also proved here.

  In \S 4, we prove miscellaneous properties of the
  geometric minimal model: functoriality,
 liftability of  group actions, the quotient property,
 the direct product property; we also relate it with Mori's minimal model
 in the 3-dimensional case.

  Throughout this paper, an analytic space $X$ is
  an analytic $\bC$-space  in Hironaka's sense \cite{hironaka},
  i.e.,

 (i) for every point $x$ of $X$,
  there exists an open neighbourhood such that the restriction of
  $X$ onto it is $\bC$-isomorphic to a local analytic $\bC$-space,

  (ii) the underlying topological space of $X$  is countable at infinity,
  i.e., \( X \) is a union of countably many compact subsets, and

  (iii) $X$  is a Hausdorff space.

Analytic spaces are always assumed to be irreducible and reduced
  except for Lemma \ref{sequence}.

  Professor Hajime Tsuji pointed out that the technique used to prove
  Theorem 3.4
  in the preliminary version is the same  as that used in Fujiki's
  paper
  \cite{fujiki2}, thus our original proof  could be  shortened by quoting
  \cite{fujiki2}, for which the authors are grateful.
    The authors also would like to thank Professors
  Kenji Matsuki, Hironobu Maeda, Noboru Nakayama and David Mathieu
  for providing  useful information.
  The authors express their
   heartfelt thanks to the referee
  for the numerous constructive comments.
%


\section{The greatest common factor of a category of  proper surjective
  morphisms}

\begin{defn}
  For an analytic space $X$,
  denote by ${\cal{PS}}(X)$ the category whose objects
  are proper surjective morphisms  $X' \to X$
  from   analytic spaces \( X' \) to \( X \)
  and whose morphisms are proper surjective morphisms over  $X$.

  In this category, we introduce an order between the objects
  according to the existence of morphisms: i.e.,
  $(X'\to X) \geq (X''\to X)$, if there is a proper surjective
  morphism $X'\to X''$ over $X$.
  In this case we say that  $(X'\to X)$ is greater than  $(X''\to X)$.
  (Here note that ``greater than'' contains  the case ``equal to.'')
\end{defn}

\begin{defn}
  For a subcategory  ${\cal F}$ of  ${\cal{PS}}(X)$,
  an object  $f:Y \to X$ of ${\cal{PS}}(X)$
  is called a {\sl common factor  } of ${\cal F}$,
  if every object of ${\cal{F}}$ is greater than   $f$.

  A maximal object among the common factors of ${\cal F}$
  is called a {\sl maximal common factor of} ${\cal F}$,
  and if a maximal common factor is unique, it is called the
  {\sl greatest common factor}.
  In a slight abuse of terminology, we sometimes say that the space
  $Y$ is a common factor or the greatest common factor
  instead of referring to the morphism $Y\to X$.

  For an object $\tilde X \to X$ and a common factor $Y\to X$ of ${\cal F}$,
  there is a morphism $\tilde X \to Y$ over $X$ by definition.
  We call this morphism a {\sl factorizing morphism}.
  If $Y\to X$  is bimeromorphic,
  then a factorizing morphism is unique.

  Note that if \( {\cal F}' \) is a subcategory of \( {\cal F} \),
  a common factor of \( {\cal F} \) is a common factor of \( {\cal F}' \).
  In particular the greatest common factor of \( {\cal F}' \) is
  greater than  the greatest common factor of \( {\cal F}
  \),
  if these greatest common factors exist.
\end{defn}

\begin{thm}
\label{existence}
  Assume that in a subcategory  ${\cal F}$ of  ${\cal{PS}}(X)$
  there is an  object $X'\to X$ which is bimeromorphic.

 Then, for every common factor \( Y'\to X \) of ${\cal F}$
there is a  greatest common factor  \( Y \to X \)  such that
 \( (Y\to X ) \geq (Y'\to X)\).

  If for every object   $X'\to X$  of  ${\cal F}$, $X'$ is normal,
  then the greatest common factor
of  ${\cal F}$ is normal.
 \end{thm}

To prove  the theorem,
  we need the following lemma and definition.

\begin{lem}
\label{sequence}
  Let $\cdots E_{i+1}\stackrel{\varphi_{i+1}}\longrightarrow E_i
   \stackrel{\varphi_{i}}\longrightarrow E_{i-1}\cdots \to E_0$
  be a sequence of proper surjective  morphisms of reduced
  (not necessarily irreducible) analytic spaces.
  For each  $i \geq 0$,
  let  $\psi_i:E \to E_i$ be a proper surjective morphism of
  reduced
  analytic spaces.
  Assume the commutativity  $\psi_i=\varphi _{i+1}\circ\psi_{i+1}$
  for every  $i\geq 0$.

  Then,
  for every  compact  subset $K\subset E_0$,
  there exists
  a number  $r\geq 0$ such that  $\varphi_i$ restricted onto
  the inverse image of
  $K$  is of
  finite fibers for all  $i\geq r$.
\end{lem}

\begin{pf}
  Since we restrict onto a compact subset,
  we may assume that the number of irreducible components of $E$ is finite.
  Then,
  we may assume that the number of irreducible components of  $E_i$'s
  is
  constant, say  $s$, because they are bounded by that of  $E$.
  Let  $E_{i}^j$  $(j=1,..,s)$
  be the irreducible components of $E_i$ and $\varphi_{i+1}(E_{i+1}^j)=E_i^j$.
  Note that for every $j$ there is an irreducible component $E^j$ of $E$ which
  is mapped onto   $E_{i}^j$.
  Then, for each $j$, we may assume that
  $\dim {E_i}^j$'s are constant for all $i$,
  because they are bounded by  $\dim E$.
  Now let us prove the lemma by  induction on $n=\dim E$.

  Assume $n=1$.
  Then   $\dim {E_i}^j=0 $ for all $i$ or  $\dim {E_i}^j=1$ for all $i$.
  Both cases are trivial.

  Let $n\geq 2$ and assume that the lemma holds true  whenever \(
  \dim E < n \).
 Since it is sufficient to prove the finiteness of $\varphi_i$ on each
  irreducible component  $E_{i}^j$',
  we may assume  that $E_i$'s and $E$ are all irreducible.
  Let  $\dim E_i=k$ and $\Phi_i=\varphi_1\cdots \varphi_{i-1}\circ\varphi_i$.
  Define $D_i\subset E_i$ as
  $\{x\in E_i \mid \dim_x \Phi_i^{-1}\circ \Phi_i(x)> 0\}$,
  then it is a proper closed subset of $E_i$.
  Next define $D\subset E$ as
  $\{y\in E \mid \dim_y\psi_0^{-1}\circ \psi_0(y)> n-k\}$,
then it is also a proper closed subset in $E$ and satisfies
  $D_i \subset \psi_i (D)$.
  This inclusion is proved as follows:

  For any point $x\in D_i$,
  take an  irreducible component $D_{i0}$ of $\Phi_i^{-1}\circ \Phi_i(x)$
  containing  $x$, then $\dim D_{i0} >0$.
  Let $D'\subset E$ be an irreducible component of \(
  \psi_{i}^{-1}(D_{i0}) \) that dominates  $D_{i0}$
  and $y\in D'$
  a point corresponding to $x$.
  Then, according to a general theory (see, for example,
  \cite[Chap.1,
  \S 8, Theorem 2]{mumford}), $\dim D' \geq \dim D_{i0}+ n-k$.
  By $D'\subset \psi_0^{-1}\circ \psi_0(y)$,
  we obtain that  $\dim_y\psi_0^{-1}\circ \psi_0(y)> n-k$,
  which yields  $x\in \psi_i(D)$.

  Now we obtain a sequence
  $\cdots \psi_i(D)\to \psi_{i-1}(D) \to \cdots\to \psi_0(D)$ of proper
surjective
  morphisms of reduced schemes which are dominated by  $D$ of dimension $< n$.
  By the induction hypothesis,
  $\psi_i(D)\to \psi_{i-1}(D)  $ is of finite fibers for $i >> 0$.
  On the other hand  $\varphi_i: E_i\to E_{i-1}$ is of finite fibers outside of
 $D_i$,
 therefore $\varphi_i$ is of finite fibers on the whole of $E_i$ for $i>>0$.
\end{pf}

\begin{defn}
  Let $K$  be a compact subset of $X$.
  We say that a common factor $Y\to X$ is {\sl maximal on } $K$,
  if for every common factor $Z\to X$ such that \( (Z\to X)\geq (Y\to
  X) \)
  we have $Z|_K=Y|_K$, where  $\ \ |_K$ means  the inverse
  image of  $K$.
\end{defn}

\noindent
{\it Proof of Theorem \ref{existence}.}
  First note that all common factors of ${\cal F}$
  are bimeromorphic to $X$ by the assumption of the theorem.
  Let $Y\to X$  be a common factor of ${\cal F}$.
  We will prove that for a compact subset $K$, there exists a common factor
  $Y'\to X$  such that \( (Y'\to X)\geq (Y\to
  X) \) and maximal on $K$.
  If $Y$ is not maximal on $K$,
  there is a common factor $Y_1\to X$ greater than $Y=Y_0$
  and not equal to $Y_0\to X$ on $K$.
  Next, if $Y_1$ is not maximal on $K$, there is $Y_2\to X$ as above.
  Successively we obtain a sequence:
 $$\cdots \to Y_i\stackrel{\varphi_i}\longrightarrow Y_{i-1}\to
  \cdots \stackrel{\varphi_1}\longrightarrow Y_0$$
   of common factors of ${\cal F}$.
  Let $X' \to X$  be an object of ${\cal F}$ and
  $\psi_i:X' \to Y_i$
   the factorizing morphism  for every $i$.
  Since $Y_i$'s are bimeromorphic to $X$,
   the commutativity  $\varphi_i\circ \psi_i=\psi_{i-1}$ follows
  for every $i$.
  Take the normalizations ${X'}^*$ and $Y_i^*$ of $X'$ and $Y_i$,
  respectively.
  Then by the lemma,
  there is an integer  $r=r(K)$ such that for every  $ i > r$
  we have
  that $Y_{i+1}^*\to Y_i^*$
  is a proper bimeromorphic morphism between normal
  spaces,
  and its restriction
  onto the inverse
  images of $K$
  is
  of finite fibers,
   which is
   therefore an isomorphism by the Zariski Main Theorem
  (for an analytic version, see \cite[4.9]{fischer}).
  Hence, on $K$, every common factor in
$$\cdots \to Y_i\stackrel{\varphi_i}\longrightarrow Y_{i-1}\to
  \cdots \stackrel{\varphi_{r+1}}\longrightarrow Y_{r}$$
 is dominated by the common normalization
  $\pi_i: Y^* \to Y_i$.
  Therefore, the coherent ${\cal O}_{Y_{r}}$-module ${\pi_r}_*{\cal O}_{Y^*}$
  has an ascending chain of submodules:
$${\cal O}_{Y_r}\subset {\rho_{r+1}}_*{\cal O}_{Y_{r+1}}\subset ...
 \subset {\rho_{i}}_*{\cal O}_{Y_{i}}\subset ...,$$
where $ \rho_{i}: Y_i \to Y_r$ is the composite  $\varphi_r\circ \cdots \circ
  \varphi_i$.
  By the coherence of ${\pi_r}_*{\cal O}_{Y^*}$
  this chain should  terminate at a finite stage.
  Thus, we obtain a common factor $Y'\to X$ greater than  $Y$ and maximal on
$K$.

  Since $X$ is countable at infinity,
  $X=\bigcup_{i=1}^{\infty}K_i$
  with $K_i$ compact and $K_i\subset K_{i+1}$.
  Take a sequence
 $$\cdots \to Y'_i\stackrel{\varphi'_i}\longrightarrow Y'_{i-1}\to
  \cdots \stackrel{\varphi'_1}\longrightarrow Y_0=Y$$
of common factors such that $Y'_i\to X$ is maximal on $K_i$.
  Then,
 $$\cdots \to Y'_i\stackrel{\varphi'_i}\longrightarrow Y'_{i-1}\to
  \cdots \stackrel{\varphi'_{n+1}}\longrightarrow Y'_n$$
are all isomorphic on $K_n$.

  Now, as in \cite[4.1.7]{flat}, we can define an analytic space
  $Y_{\infty}=\underset{r}{\underset{\longleftarrow}{\lim}}Y'_r$,
  which is again a common factor of ${\cal F}$  and greater than all
  $Y'_i$'s.
 This $Y_{\infty}$ is a maximal common factor of ${\cal F}$.
  Indeed, if $Z\to X$ is a common factor greater than $Y_{\infty}$,
  then for every $K_i$,
  $Z|_{K_i}=Y'_i|_{K_i}=Y_{\infty}|_{K_i}$ by the maximality of $Y'_i$ on
$K_i$.
  Hence $Z=Y_{\infty}$.

  Let  $Y\to X$ and $Y'\to X$ be  maximal common factors.
  Then the unique component $Z$ of $Y\times_X Y'$ dominating $Y$ and $Y'$
  is again a common factor of ${\cal F}$ and greater than
  $Y$ and $Y'$,
  because of the universality of the fiber product.
  By the maximality of $Y$ and $Y'$, $Z$ coincides with $Y$  and $Y'$,
  which yields the uniqueness of
  the greatest common factor.

  If every object of ${\cal F}$ is normal, the normalization of a common factor
  is also a common factor.
  So the greatest common factor is normal.
$\Box$


\begin{rem}
\label{hilb}
Instead of Lemma \ref{sequence}, we can prove the theorem  differently:
Consider a sequence of common factors as in the proof of Theorem
\ref{existence}.
Each $Y_i$ is dominated by a  common object
$\tilde X\stackrel{f_i}\longrightarrow Y_i$.
 Let ${\cal J}_{i,a}$ be the defining ideal of
  ${f_i^{-1}f_i(a)}$ in $\tilde X$.
 On a compact subset $K$ on $Y_0$, consider the Hilbert-Samuel function
  $H_{i,a}$ of ${\cal O}_{\tilde X, a}/{\cal J}_{i,a}$ for $a\in f_0^{-1}(K)$.
   One can  prove that there exists $n(K)$ such that for any
   $n\geq n(K)$ and $a\in f_0^{-1}(K)$, the equality $H_{n,a}=H_{n(K), a}$
holds
  by means of Zariski upper-semicontinuity (\cite[lemmas 7.1, 7.2]{bm1},
\cite[Lemma 3.10 and
Definition 3.11]{b-m})
  and the stabilization property
  of Hilbert-Samuel function (\cite[Theorem 5.2.1]{bm4}).
  This equality yields the equality of the ideals ${\cal J}_{n, a}={\cal
J}_{n(K),a}$,
  and therefore $Y_n\to Y_{n(K)}$ is finite on $K$.
  The rest of the proof follows in the same way as  above.
 \end{rem}

\begin{rem}
\label{alg}
  Consider the case that $X$ is an algebraic variety.
  Theorem \ref{existence} holds also for the category ${\cal{PS}}alg(X)$
  whose objects are  proper surjective morphisms of algebraic
  varieties to $X$.
  The proof is easier.
  In Lemma \ref{sequence},
  there is a number $r$ so that $\varphi_i$ $(i\geq r)$ is of finite fibers
  for the global  $E_i$'s.
  Therefore one can prove Theorem \ref{existence} immediately by showing the
  ascending chain condition,  without using $Y_{\infty}$.
\end{rem}

\begin{exmp}
  There are many subcategories of ${\cal{PS}}(X) $  which satisfy the
condition of the theorem.

  ${\mathcal R}es(X)$ : the category of resolutions of the singularities of
   $X$.
   Here a resolution \( f:\tilde X \to X \) means that \( f \) is a
   proper bimeromorphic morphism from non-singular analytic space
   \( \tilde X \).

  ${\mathcal S}m(X)$ : the category of proper surjective morphisms
  $X'\to X$,
  where $X'$ is non-singular.
  Here, we note that the dimension of $X'$ may be bigger than the dimension
  of $X$.

  ${\mathcal C}(X)$  (resp. ${\mathcal L}t(X)$, ${\mathcal T}(X)$) :
   the category of proper bimeromorphic  morphisms
    $X'\to X$ where \( X' \) has  canonical singularities (resp. log-terminal
   singularities, terminal singularities).

  These are of course nonempty categories due to the results of
  \cite{hironaka}, \cite{hironaka2},  \cite{bm2}, \cite{b-m}  and
  \cite{vil}.

  We can also think of the following categories for a toric variety $X$
(cf., \ref{alg}):

${\mathcal R}es(X)_{tor}$ : the category of toric resolutions of the
singularities of $X$.

 ${\mathcal S}m(X)_{tor}$ : the category of toric proper surjective
morphisms  $X'\to X$,
  where $X'$ is non-singular.

  We do not know if the greatest common factors of  ${\mathcal S}m(X)$ and
of ${\mathcal S}m(X)_{tor}$ are different  for the same $X$.

\end{exmp}

\begin{exmp}
  For a two-dimenstional analytic space $X$,
  the greatest common factor of ${\mathcal R}es(X)$ is the minimal resolution.
\end{exmp}

\begin{defn}
\label{def of gmm}
  The greatest common factor of  ${\mathcal S}m(X)$ is called the {\sl
geometric
minimal
   model} of $X$.
\end{defn}

\begin{defn}
  For an algebraic variety $X$,
  the greatest common factor of ${\cal S}m(X)$ in the category  ${\cal
PS}alg(X)$ is
  called the {\sl algebraic geometric minimal model} which is not studied in
the
  following
  sections (cf.  \ref{alg}).
\end{defn}


\section{The geometric minimal model and rational curves}

\noindent
  In this section we study some properties of the geometric minimal model of
an analytic
 space    $X$.

  First, we state and prove  a basic lemma,
  which is well known and used often in this paper.

\begin{lem}
\label{key}
  Let $f:Z\to X$ be a proper surjective morphism and
  $g:Z \to Y$ a  morphism.
  Assume  $X$  is normal and
 $g$  is constant on $f^{-1}(x)$ for every $x\in X$.

  Then, there is a morphism $\varphi: X \to Y$
  such that $g=\varphi\circ f$.

\end{lem}

\begin{pf}
  Let $\Gamma\subset X\times Y$ be the reduced image of the proper morphism
  $(f,g): Z \to X\times Y$ and $\Phi:\Gamma \to X$
  the restriction of the projection
  $X\times Y \to X$  onto $\Gamma$.
  Then,  $\Phi^{-1}(x)$
   consists of one point.
  Now we obtain that $\Phi:\Gamma \to X$ is a bijective
  proper bimeromorphic morphism.
   Since
  $X$ is normal,  $\Phi$ is an isomorphism by  Zariski's Main Theorem.
\end{pf}

The following is a basic property of the
geometric minimal model.

\begin{thm}
\label{main thm}
  Let $X$ be an analytic   space of dimension $\geq 1$,
  $\varphi: Y\to X$ a common factor of ${\mathcal S}m(X)$
  and $f:X'\to X$
  an arbitrary object of ${\mathcal S}m(X)$.

  Then every rational curve $l\subset X'$ which is mapped to a point on $X$ is
  mapped to a point on $Y$ by the factorizing morphism $X'\to Y$.
\end{thm}

\begin{pf}
  It is sufficient to construct a commutative diagram
  with $f$, $\gamma$ and $f\circ \alpha=\gamma\circ \beta$ being
in ${\mathcal S}m(X)$
$$
\begin{matrix}
\tilde X& \stackrel{\beta}\longrightarrow & X''\\
\alpha \downarrow &  &  \gamma\downarrow\\
X' &  \stackrel{f}\longrightarrow & X\\
\end{matrix}
$$
  such that there is a curve $\tilde l$ on $\tilde X$ with $\alpha(\tilde l)=l$
  and $\beta(\tilde l)=p$ a point.
  This is because $\gamma $ and $f$ factor through $\varphi:Y\to X$, and
  the image of $l$ on $Y$ coincides with the image of $\tilde l$ and therefore
  with the image of $p$.

 Now we construct $\tilde X$ and $\tilde l$.
  First, if $\dim X'=1$, then there is no such rational curve and
  therefore the statement is trivial.
  If $\dim X'=2$, take the product $X_{0}=X'\times \bP^1$.
  Then, the composite  $X_{0}=X'\times \bP^1\stackrel{p_1}\longrightarrow X'
\to X$
  is also an object of ${\cal  S}m(X)$,
  where $p_1$ is the projection to the first factor.
  Take a section $l_{0}$ of $p_1$ over $l$.
  If \( \dim X'\geq 3 \), let \( X_{0}= X' \) and \( l_{0}=l \).

  Then for \( \dim X'\geq 2 \),
  there is a proper bimeromorphic morphism  \(\alpha_{1}: X_{1}\to
  X_{0} \) with
  a non-singular rational curve \( l_{1}\subset X_{1} \) such that
  \( \alpha_{1}(l_{1})=
  l_{0} \) and the normal sheaf \( N_{l_{1}/X_{1}}={\cal O}_
  {{\bP}^1}(-1)\oplus {\cal O}_{{\bP}^1}(-1)\oplus \cdots \oplus
  {\cal O}_{{\bP}^1}(-1)  \) by Fujiki's result \cite[Lemma 5]{fujiki2}.
  The following construction is also found in \cite{fujiki2}:
  take the blow-up \(\alpha_{2}: \tilde X \to X_{1} \) with the center
  \( l_{1} \) and define \( \alpha:=\alpha_{1}\circ \alpha_{2} \).
      Then the exceptional divisor $E$ of $\alpha_{2}$ is
$\bP^1\times\bP^{n-2}$
and its normal sheaf
   is $p_1^*{\cal O}_{\bP^1}(-1) \otimes p_2^*{\cal O}_{\bP^{n-2}}(-1)$,
  where $p_i$ is the projection to the $i$-th factor and \( n \) is
  the dimension of \( X' \).
  Then,  there is a proper bimeromorphic morphism
  $\beta:\tilde X \to X''$
  whose restriction onto $E$ is the projection
  $E \stackrel{p_2} \longrightarrow \bP^{n-2}$.
  By the normal bundle of $E$, we can see that $X''$ is non-singular
   (see, for example, \cite[Theorem 3, (2)]{ishii}).
  Since $f\circ \alpha$ is constant on $\beta^{-1}(x)$ for every
  $x\in X''$,
  we obtain a morphism
  $\gamma:X''\to X$ with $\gamma\circ \beta=f\circ \alpha$ by Lemma \ref{key}.
  Now,   take a section $\tilde l=\bP^1\times \{P\}\subset E$ over $l_1$,
  where $P$ is a point of $\bP^{n-2}$.
 Then $\beta(\tilde l)$ is a point, which shows that $\tilde l$ is a
 curve  we wanted.

\end{pf}

\begin{thm}
\label{res}
  Let $X$ be an analytic   space of dimension $\geq 3$,
  $\varphi': Y'\to X$ a common factor of ${\cal R}es(X)$
  and $f:X'\to X$
  an arbitrary object of ${\mathcal R}es(X)$.

  Then every rational curve $l\subset X'$ which is mapped to a point on $X$ is
  mapped to a point on $Y'$ by the factorizing morphism $X'\to Y'$.
\end{thm}

\begin{pf}
  Since all objects in ${\cal R}es(X)$ are of dimension $\geq 3$,
  we can construct a commutative diagram in ${\cal R}es(X)$
  in the same way as
  in the proof of the
  previous theorem.
\end{pf}

\begin{thm}
\label{property}
    Let $Y \to X$ be a bimeromorphic object of ${\cal{PS}}(X)$ satisfying the
   following
   condition (R):

  (R)
  There exists a proper bimeromorphic morphism $\tilde Y \to Y$  such that
    every rational curve $l \subset \tilde Y$ which is mapped to a point in $X$
    is mapped to a point in $Y$.

   Then $Y \to X$ is a common factor of ${\cal  S}m(X)$.
\end{thm}

\begin{pf}
  Let $X'\to X$ be an arbitrary object of ${\cal  S}m(X)$.
 Then, there is a component $Z$ of $X'\times_X \tilde Y$ on which
 $X'\times_X \tilde Y \to X'$ is bimeromorphic.
  Then, by \cite[Main Theorem II']{hironaka} and \cite[Corollary 1]{flat},
  for every point $x\in X'$ there are a relatively compact  neighbourhood
  $U_r$, $(r\geq 1)$ of $x$, a non-singular
  analytic space $ U_1$, bimeromorphic morphisms
  $f: U_1 \to U_r$ and $ U_1 \to Z|_{U_r} $
  such that $f$  is the composite $f_{r-1}\circ f_{r-2}\circ\cdots \circ f_1$
  of blow-ups $f_i: U_i\to  U_{i+1}$
  with non-singular centers and the following diagram is commutative:
$$
\begin{matrix}
 U_1     & \stackrel{f_1}\longrightarrow &  U_2 & \to &
      ...&\stackrel{f_{r-1}}\longrightarrow  &  U_{r} & \subset & X'\\
\downarrow &                             &            &
         &                               &            &     &   \\
Z        &                               &             &  &
         &                               &            &     & \downarrow\\
\downarrow &                             &            &     &
          &                              &            &     &    \\
\tilde Y  &  \to                         & Y          &     &
          & \to                          &            &     & X.
\end{matrix}
$$

  Therefore,  $f_1^{-1}(y) = \bP^s$  $(s\geq 0)$
  for every $y\in U_2$.

  For a rational curve $l \subset f_1^{-1}(y)\subset  U_1$,
  if the image of $l$ on $\tilde Y$ is not a point then it is a rational curve.
  So  $l$  is mapped to a point in $Y$ by the condition (R), which yields
  the existence of a morphism $U_2\to Y$ by Lemma \ref{key}.
  In the same way we induce the  morphisms $ U_i\to Y$ for $i=2,.., r$.
  Now we complete the proof of the statement that $ X' \to X$ factors
through $Y\to X$.
\end{pf}

\begin{cor}
\label{coincide}
  Let $X$ be an analytic   space of dimension $\neq 2$.

  Then, the greatest common factor of ${\cal R}es(X)$
  coincides with the geometric minimal model.
\end{cor}

\begin{pf}
   Since  ${\cal R}es(X)\subset {\cal  S}m(X)$, it follows that
  the greatest common factor $Y'\to X$ of ${\cal R}es(X)$ is greater than
the geometric minimal model
  $Y\to X$.
  For $\dim X=1$, $Y$ and $Y'$ are normal and bimeromorphic to $X$,
  therefore they coincide by Theorem \ref{existence}.

  For $\dim X\geq 3$, it is sufficient to prove that $Y'\to X$ is a common
factor
  of ${\cal  S}m(X)$ by Theorem \ref{existence}.
  This follows from Theorem \ref{res} and Theorem \ref{property}.
\end{pf}

\begin{thm}
   The geometric minimal model of an analytic   space $X$
 is the maximal bimeromorphic object among
  the objects satisfying $(R)$  in ${\cal{PS}}(X)$.
\end{thm}

\begin{pf}
  This is clear by Theorem \ref{existence} and Theorem \ref{property},
  because the geometric minimal model also satisfies $(R)$
  by Theorem \ref{main thm}.
\end{pf}

\begin{cor}
     If $Y\to X$ is the geometric minimal model,
   then    for every non-trivial bimeromorphic morphism $Z\to Y$,
  there is a rational curve on $Z$ which is mapped  to a point on $Y$.
\end{cor}

\begin{pf}
  If $Y\to X$ is the geometric minimal model,
  then $Z$ is not a common factor of ${\cal S}m(X)$ by the maximality of $Y$.
  Therefore, by Theorem \ref{property},
  for a resolution $\tilde Z \to Z$ of singularities of $Z$
  there is a rational curve $l \subset \tilde Z$ which is mapped to a point
on $X$
  and to a curve $l'$ on $Z$.
  This curve $l'$ is a  rational curve on $Z$, as required.
\end{pf}

  Now, the above theorem and Theorem \ref{existence} imply the
  following:

\begin{cor}
  If $X$ is an analytic   space of dimension two,
  then the geometric minimal model of $X$ is a surface obtained
  by contracting all rational exceptional curves of the minimal
  resolution.
\end{cor}

\begin{exmp}
  Let $X$ be an affine cone over a non-singular variety $A$,
  which means that there is a  line bundle $Y$ over $A$
  with the zero section of negative normal bundle and the zero
  section is contracted to a point of \( X \) in \( Y \).
  We will denote the zero section of \( Y \) also by \( A \).

  (1) Suppose \( A \) is either an Abelian variety or the direct product
  of
  non-singular curves of positive genera.
  Then $Y$ is a resolution of the singularities on $X$ and
  also the geometric minimal model of \( X \), since \( A \) has no
  rational curve (see Theorem \ref{property}).

  (2) If \( A \) is an algebraic K3-surface, then there is a rational curve
  with  positive self-intersection number.
  Then, the whole of \( A \) should be contracted to a point in the
  geometric minimal model; therefore the geometric minimal model of \( X \)
  is the identity.

  (3) Let \( A \) be a ruled surface with \( q\geq 1 \) and
  \( p:A \to C \)  the ruling,
  then there is a contraction
  \( Y \to Y' \)  of the zero section \( A \) to a curve \( C \) in \( Y \)
  by  \cite[Theorem 1]{fujiki}.
  Then \( Y' \) is the geometric minimal model of \( X \).

\end{exmp}

\begin{exmp}
  Let \( X \subset \bC^{4}\) be a hypersurface defined by a homogeneous
  polynomial \( f \) such that the hyperplane in \( \bP^{3} \) defined
  by \( f \) is non-singular and hyperbolic (for the definition of
  hyperbolic, see  \S 4 of \cite{fujiki2}).
  Then, the blow-up \(\varphi: Y\to X \) of \( X \) at the origin is the
  geometric minimal model of \( X \), because there is no rational
  curve in a fiber of \( \varphi \) and \( Y \) is non-singular.
\end{exmp}

\begin{exmp}
  Let $X$ be a toric variety,
  then the geometric minimal model is the identity.
  Indeed, for every toric resolution $f:X' \to X$,
  each fiber $f^{-1}(x)$ of $f$ is a connected union of rational varieties,
  so it should be mapped to a point on the geometric minimal model.
\end{exmp}



\section{Miscellaneous properties on the geometric minimal models}

\begin{prop}(Functoriality)
\label{functoriality}
  Let $X$ and $X'$ be   analytic  spaces
  and $\varphi:Y_X\to X$ and $\varphi':Y_{X'}\to X'$
  the geometric minimal models of
  $X$ and $X'$, respectively.

  Then a proper surjective morphism
  $f:X\to X'$ can be uniquely lifted to $\tilde f:Y_X\to Y_{X'}$ such that
  $\varphi'\circ \tilde f=f\circ \varphi$.
\end{prop}
\begin{pf}
  The composition with \( f \) gives an injection \( {\cal S}m(X)\to
  {\cal S}m(X') \).
  Therefore the greatest common factor of the image of \( {\cal S}m(X)
  \) is greater than that of \( {\cal S}m(X') \).
\end{pf}

\begin{cor}
  Let $\pi:X\to S$ be a proper surjective morphism to a surface $S$.
  Let $S'$  be a surface obtained by contracting all rational exceptional
  curves of the minimal resolution of $S$.
  If the geometric minimal model of $X$ is the identity
  (for example, $X$ is toric),
  then $\pi$ factors through $S'$.
\end{cor}

\begin{cor}
   Let $X$ be  an analytic   space and
   $\varphi:Y\to X$  the geometric minimal model of $X$.

   Then every automorphism $\sigma: X \to X$ can be lifted to $\sigma':Y\to Y$
  such that $\varphi\circ \sigma'=\sigma\circ \varphi$.
\end{cor}

\begin{pf}
  By Proposition \ref{functoriality},
  $\sigma$ can be lifted to $\sigma ':Y\to Y$ and $\sigma^{-1}$ can be
lifted to
  $\sigma'':Y\to Y$.
  By the uniqueness of the lifting, $\sigma'$ and $\sigma''$ are the inverse
  morphism of each other.
\end{pf}

  By this corollary, we obtain that every action of a finite group on $X$ is
  lifted to the geometric minimal model.
  What follows is a discussion of the quotient by a finite group.

\begin{prop}
\label{quotient}
   Let $X$ be
   an analytic   space and $G$ a finite group acting on $X$.
   Let $\varphi: Y_X\to X$ and $\varphi':Y_{X/G}\to X/G$ be the geometric
   minimal models of $X$ and $X/G$, respectively.

 Then there is a morphism $Y_X/G \to Y_{X/G}$ over $X/G$.
\end{prop}

\begin{pf}
  By the functoriality of the geometric minimal models (cf. Proposition
\ref{functoriality}), there is a morphism $Y_X\to Y_{X/G}$.
  By the uniqueness of this morphism, it is a $G$-equivariant morphism
 with the trivial action on $Y_{X/G}$.
  Hence the morphism factors through $Y_X/G$.
\end{pf}

  The following is obvious from the  above proposition.

\begin{cor}
\label{quot}
  Let $X$ be  an analytic   space and $G$ a finite group acting on $X$.
  If the geometric minimal model of $X$ is the identity,
  then the geometric minimal model of $X/G$ is also the identity.
\end{cor}
\begin{exmp}
  Let $X=\bC^n/G$ by a finite group $G$,
  then the geometric minimal model of $X$ is the identity.
\end{exmp}

\begin{exmp}
   The equality $Y_X/G \simeq Y_{X/G}$ does not hold in general.
  Indeed, take a line bundle $Y'={\bf V}({\cal O}_{\bP^1}(m))$ over $\bP^1$,
  where $m$ is a positive integer.
  Let $\pi:C \to \bP^1$ be a double cover from an elliptic curve $C$
  and $G=\bZ/(2)$ the  Galois group.
  Then $G$ acts on $Y={\bf V}(\pi^*{\cal O}_{\bP^1}(m)  )$.
  Denoting by  $X$ the variety obtained by contracting the zero section
  of $Y$, the group \( G \) acts on \( X \) as well.
  The quotient  $X/G$ has a resolution $Y'$, therefore it is the quotient
of $\bC^2$
  by a finite group of order $m$.
 Hence the geometric minimal model of $X/G$ is the identity,
  but the geometric minimal model of $X$ is $Y$, and $Y/G\simeq Y'$.
\end{exmp}

  In order to prove the direct product property of the geometric minimal model,
  we state the following lemma.

\begin{lem}
\label{product}
  Let $X$ and $X'$ be analytic  spaces,
  ${\cal F}$ a subcategory of  ${\cal{PS}}(X)$
  satisfying the condition of Theorem \ref{existence}, i.e., there is
  a bimeromorphic object \( X' \to X \) in \( {\cal F} \).
  Let
  $Y_X$ be the greatest common factor of ${\cal F}$ and
  $f:Z\to X'$ an object of ${\cal{PS}}(X')$.
  Let ${\cal F}\times Z$ be the category of the morphisms
  $\tilde X\times Z\to X\times X'$  with $\tilde X\in {\cal F}$.

   Then $Y_X\times Z$ is the greatest common factor of ${\cal F}\times Z$.
\end{lem}

\begin{pf}
  It is clear that $Y_X\times Z$ is a common factor of
  ${\cal F}\times Z$.
  By Theorem \ref{existence},
  there is the greatest common factor $Y'$ of
  ${\cal F}\times Z$.
  So for every object $\tilde X\to X$ of ${\cal F}$,
  there is a factorizing morphism $\varphi:\tilde X\times Z \to Y'$
  of $\tilde X \times Z \to X \times X'$.
  On the other hand, there is a morphism $\psi:Y'\to Y_X\times Z$
  by the maximality of $Y'$.
  Therefore for every point $z\in Z$,
  $\tilde X\times\{z\}\to X\times \{f(z)\}$ is factored as follows:

$$\tilde X\times\{z\} \stackrel{\varphi}
\longrightarrow
  \psi^{-1}(Y_X\times\{z\})\stackrel{\psi }\longrightarrow
  Y_X\times \{z\} \to
 X\times \{f(z)\}.$$

  So, \( \psi^{-1}(Y_X\times\{z\}) \) is a common factor of \( {\cal
  F} \).
  Then, by the maximality of $Y_X$,
  $\psi^{-1}(Y_X\times\{z\})\stackrel{\psi }\longrightarrow
  Y_X\times \{z\}$ is an isomorphism.
  Now we obtain that $\psi:Y'\to Y_X\times Z$ is bijective and
  $\psi^{-1}((y, z))\simeq \{(y,z)\}$ for every point $(y,z)\in Y_X\times Z$,
  which yields that  $\psi$ is an isomorphism.
\end{pf}

\begin{prop}(Direct product)
\label{direct product}
  Let $X$ and $X'$ be analytic  spaces,
  $Y_X$, $Y_{X'}$ and $Y_{X\times X'}$ the geometric minimal models
  of $X$, $X'$ and $X\times X'$, respectively.

  Then $Y_X\times Y_{X'}\simeq Y_{X\times X'}$ over $X\times X'$.
\end{prop}

\begin{pf}
  First for any objects $\tilde X\to X$ and $\tilde X'\to X'$ of
  ${\cal S}m(X)$ and ${\cal S}m(X')$, respectively,
  $\tilde X\times \tilde X'\to X\times X'$ factors through
  $Y_{X\times X'}$.
  For every $\tilde X'$
  in ${\cal S}m(X')$,
  $Y_X\times \tilde X'$ is the greatest common factor
  of ${\cal S}m(X)\times \tilde X'$  by Lemma \ref{product}.
  So there is a morphism $Y_X\times \tilde X'\to Y_{X\times X'}$.
  Then  by Lemma \ref{product} again,
  $Y_X\times Y_{X'}$ is the maximal
  common factor of $Y_X\times {\cal S}m(X')$.
  So there is a morphism $Y_X\times Y_{X'}\to  Y_{X\times X'}$.

  In order to prove the existence of the inverse morphism $  Y_{X\times X'}\to
Y_X\times
Y_{X'}  $,
  it is sufficient to prove that $Y_X\times Y_{X'}$ is a common factor of
  ${\cal S}m(X\times X')$.
  Let $\tilde X$ and $\tilde X'$ be  resolutions of $Y_X$ and $Y_{X'}$,
  respectively, and $l$ a rational curve on
  $\tilde X\times \tilde X'$ which is mapped to a point on $X\times X'$.
  By Theorem \ref{property}, it is sufficient to prove that $l$ is mapped to
  a point on $Y_X\times Y_{X'}$.
  If the image of $ l$ by the projection
 $ \tilde X\times \tilde
   X'\stackrel{p_1}\longrightarrow \tilde X$
  is not a point, then it is a rational curve, and therefore is mapped to a
  point on $Y_X$.
  In the same way, $l$  is mapped to a point on $Y_{X'}$ through the
projection
  $\tilde X\times \tilde
  X'\stackrel{p_2}\longrightarrow \tilde X'$.
    Hence $ l$ is mapped to a point on $Y_X\times Y_{X'}$.
\end{pf}


  Next we consider an action of a complex Lie group.

\begin{prop}
  Let $G$ be a complex Lie group  acting on an
analytic   space $X$.

  Then the action can be lifted to the geometric minimal model of $X$.
\end{prop}

\begin{pf}
  Let $\sigma:G\times X \to X$ be the morphism defining the action,
  \( p_{1}:G\times X \to G \) the projection to the first factor
  and $Y$  the geometric minimal model of $X$.
  Since the geometric minimal model of $G\times X$  is $G\times Y$
  by Proposition \ref{direct product}, the automorphism
  $(p_{1},\sigma): G\times X \to G\times X$ can be lifted to
  $\Phi :G\times Y \tilde\to G\times Y$.
  The composite \( p_{2}\circ\Phi:G\times Y \to Y \)
   defines the lifted action of $G$ on $Y$,
   where \( p_{2}:G\times Y \to Y \) is the projection to the second
   factor.
\end{pf}

\begin{say}
  It is natural to ask if the geometric minimal model is a local property,
 i.e., if $\varphi: Y\to X$ is the geometric minimal model of $X$ and
$U\subset X$
 is open, then $\varphi^{-1}(U)\to U$ is the geometric minimal model of $U$.
  The following example shows that the geometric minimal model is not local.
\end{say}

\begin{exmp}
  Here we construct  $U\subset X$ such that the geometric minimal model of $X$
  is the identity but that of $U$ is not.

Let $S$ be a compact rational surface with only one simple elliptic
singularity $p$ on it.
  Let $V=S\times \bP^1$ and $C\not\ni p$ be a non-singular very ample
divisor on $S_0=S\times\{0\}$.
  Take the blow-up $X'\to V$ with the center $C$,
  then the proper transform $S'_0$ of $S_0$ on $X'$ has the negative normal
bundle on $X'$.
  Hence, there is a  morphism  $f:X'\to X$  which contracts $S'_0$ to a
point $Q$ of
   $X$.
  Since $X'$ has singularities
  (simple elliptic singularities)$\times \bP^1$, there is a resolution
  $g:\tilde X\to X'$ whose exceptional divisor $E$ is isomorphic to
(elliptic curve)$\times \bP^1$.
  Then, the composite $f\circ g:\tilde X \to X$ factors through the
geometric minimal model
  $Y\to X$ of $X$.
  The proper transform $S''_0$ of $S'_0$ on $\tilde X$ is contracted to a
point on $Y$,
  because $S''_0$ is rational.
  As $S''_0\cap E$ has the trivial normal bundle in $E$, the other elliptic
fibers are all contracted
  to points on $Y$, therefore $Y$ must be isomorphic to $X$.

  On the other hand, for $U=X\setminus \{Q\}$, the restriction $(f\circ
g)^{-1}(U)\to U$ is a
  resolution of singularities on $U$ and has no rational curves in each fiber
  of \( f\circ g \),
  therefore it is the geometric minimal model.
\end{exmp}

\begin{prop}
\label{canonical}
  Let $X$ be a three-dimensional analytic space with at worst
  canonical singularities.

  Then the geometric minimal model of $X$ is the identity.
\end{prop}

\begin{pf}
  By \cite[Corollary 2.14]{reid}, every exceptional divisor of
  a resolution $f:\tilde X \to X$
  is a ruled surface.
  Indeed, \cite[Corollary 2.14]{reid} shows that there exists a
  resolution \( X'_{1}\to X' \) of a canonical singularity \( (X',x) \)
   of index one such that all
  exceptional divisors are ruled.
  For an arbitrary resolution \( X'_{2}\to X' \),
  there exists a resolution
  \( \tilde X'\to X' \) dominating both \( X'_{1} \) and \( X'_{2} \)
  such that the morphism \( \tilde X' \to X'_{1} \) is the composite
  of blow-ups at non-singular centers.
  Therefore, all exceptional divisors for \( \tilde X'\to X' \) are ruled
  which yields that all exceptional divisors for \( X'_{2}\to X' \) are
  ruled.
  For a canonical singularity of index \( > 1 \),
  take the canonical cover and reduce the problem to the case of
  index one.
  Here, we should note that a surface which is a surjective image of a
  ruled surface is again ruled.
  As all exceptional divisors are ruled,
  the geometric minimal model $\varphi:Y\to X$ of $X$ is either   an
isomorphism or
  a contraction of the union $C$ of finitely many curves.
  Assume that $\varphi$ is not an isomorphism and let $g:\tilde X \to Y$ be
the
factorizing morphism.
  Then, all components of $C$ are rational,
otherwise
 $R^1\varphi_*{\cal O}_Y\neq 0$,
 since there is a surjection \( R^1\varphi_*{\cal O}_Y\to H^1(C_{i},
 {\cal O}_{C_{i}}) \), where \( C_{i} \) is a component of \( C \).
  This non-vanishing contradicts the edge sequence
  \[ 0\to E^{1,0}_{2}\to E^1 \to \ldots \]
  of the Leray spectral sequence:
  $$E^{p,q}_2=R^p\varphi_*R^qg_*{\cal O}_{\tilde X} \Longrightarrow
E^{p+q}=R^{p+q}f_*{\cal O}_{\tilde  X}.$$
  If the curve $C$ is not contained in the singular locus of \( Y \),
  then we can construct a resolution $\tilde X'\to X$ with a rational curve $l$
  on $\tilde X'$ mapped to a point on $X$ and mapped to a curve on $Y$, a
contradiction to Theorem \ref{main thm}.
  Thus, $C$ is contained in the singular locus of $Y$.
  Note that $Y$  has at worst canonical singularities,
  since \( K_{Y}=\varphi^*K_{X} \)
  which follows from that
  \( \varphi: Y\to X \) is isomorphic in codimension one.
  So general points of $C$ are cDV singularities in $Y$.
  For an irreducible component $C_i$ of $C$,
  take the blow-up $h:Y'\to Y$ with the center $C_i$.

  Then,
   in case $Y$ has compound $A_1$, $D_n$ $(n\geq 4)$, $E_6$, $E_7$, $E_8$
  singularities on general points of $C_i$,
 $h$ has an irreducible exceptional divisor $E$ and $h^{-1}(p)$ is an
irreducible
  rational curve for a general $p\in C_i$.
  This is because the restriction of \( h \) onto  a general hyperplane cut
at a general point of
   \( C_{i} \)
  is the blow-up of the surface at the point.
  Then, $E$ is rational and is not contained in the singular locus of $Y'$.
  Therefore, by Theorem \ref{main thm}, $E$ must be mapped to a point on $Y$,
  a contradiction.

  In case $Y$ has compound $A_n$ $(n\geq 2)$ singularities at general
points of $C_i$,
  then $h^{-1}(p)$ is the union of two $\bP^1$ with one intersection point.
  This implies that the exceptional divisor $E$ of $h$ has a generic section
over $C_i$.
Then, take the blow-up \( h_{2}:Y_{2}\to Y' \)  at the intersection curve
  of \( E \) (this curve is the generic section over \( C_{i} \)).
   Then, again blow-up \( Y_{2} \) at the intersection curve of
   the new exceptional divisor, if the intersection curve  exists.
  Continuing this procedure,
   we obtain a sequence of \( ([n/2]+1) \)-blow-ups at generic
  sections over
  \( C_{i} \):
  \[ Y_{[n/2]+1}\stackrel{h_{[n/2]+1}}\longrightarrow Y_{[n/2]}\to
  \ldots\stackrel{h_{2}}\longrightarrow Y'
  \stackrel{h}\longrightarrow Y. \]
  We can see that the exceptional divisor \( E_{[n/2]+1} \)
   of the last blow-up
  is irreducible and the fiber \( {h_{[n/2]+1}}^{-1}\circ \cdots
  \circ h^{-1}(p) \) on \( E_{[n/2]+1} \) of a general point
  \( p\in C_{i} \) is \( \bP^1 \).
  It follows that $E_{[n/2]+1}$ is rational and is not contained in the
singular locus
of $Y_{[n/2]+1}$,
  a contradiction in the same way as above.
\end{pf}

\begin{cor}
\label{log-terminal}
  Let $X$ be a sufficiently small neighbourhood of
  a  log-terminal singularity of dimension three.

  Then the geometric minimal model of $X$ is the identity.
\end{cor}

\begin{pf}
  Taking the canonical cover,
  we may assume that \( X \) is a cyclic quotient of an analytic space
  with at worst canonical singularities.
  Then, by Corollary \ref{quot} and Proposition \ref{canonical} we obtain
the statement.
\end{pf}

\begin{cor}
  Let $X$ be an analytic space of dimension three,
  $Y'\to X$  the relative canonical model in  Mori's sense
  and $Y\to X$ the geometric minimal model.

  Then there is a morphism $\varphi:Y'\to Y$ over $X$.
  Therefore, there is a morphism \( Y''\to Y \) over \( X \),
  where \( Y'' \) is a relative minimal model in Mori's sense.
\end{cor}

\begin{pf}
By \ref{functoriality},
  the geometric minimal model of \( Y' \) is greater than \( Y \) over \( X \).
  Here, by Proposition \ref{canonical}, the geometric minimal model
  of \( Y' \) is the identity.

\end{pf}

\makeatletter
\makeatother

\end{document}